\documentclass{elsart}
\usepackage{amssymb}

\usepackage[english,francais]{babel}
\newcommand{\showcomments}{yes}

\newsavebox{\commentbox}
{\ifthenelse{\equal{\showcomments}{yes}}%
{\footnotemark
        \begin{lrbox}{\commentbox}
        \begin{minipage}[t]{1.25in}\raggedright\sffamily\tiny
        \footnotemark[\arabic{footnote}]}
{\begin{lrbox}{\commentbox}}}%
{\ifthenelse{\equal{\showcomments}{yes}}%
{\end{minipage}\end{lrbox}\marginpar{\usebox{\commentbox}}}
{\end{lrbox}}}

\newcommand{\eps}{\varepsilon}

\newcommand{\Z}{\mathbb{Z}}

\DeclareFontFamily{OT1}{rsfs}{}

\DeclareFontShape{OT1}{rsfs}{n}{it}{<-> rsfs10}{}
\DeclareMathAlphabet{\mathscr}{OT1}{rsfs}{n}{it}

\newcommand{\R}{\mathbb{R}}

\newtheorem{theorem}{Theorem}[section]
\newtheorem{lemma}[theorem]{Lemma}
\newtheorem{e-proposition}[theorem]{Proposition}
\newtheorem{corollary}[theorem]{Corollary}
\newtheorem{e-definition}[theorem]{Definition\rm}

\newtheorem{theoreme}{Th\'eor\`eme}[section]

\newcommand{\vol}{\mathrm{vol}}

\newcommand {\comment} [1] {}

\def\og{\leavevmode\raise.3ex\hbox{$\scriptscriptstyle\langle\!\langle$~}}
\def\fg{\leavevmode\raise.3ex\hbox{~$\!\scriptscriptstyle\,\rangle\!\rangle$}}

\journal{the Acad\'emie des sciences}
\begin{document}

\centerline{}
\begin{frontmatter}

\selectlanguage{english}
\title{On the growth of Betti numbers of locally symmetric spaces}

\selectlanguage{english}
\author[abert]{Miklos Abert},
\ead{karinthy@gmail.com}
\author[bergeron]{Nicolas Bergeron},
\ead{bergeron@math.jussieu.fr}
\author[biringer]{Ian Biringer},
\ead{ianbiringer@gmail.com}
\author[gelander]{Tsachik Gelander},
\ead{tsachik.gelander@gmail.com}
\author[nikolov]{Nikolay Nikolov},
\ead{zarkuon@gmail.com}
\author[bergeron]{Jean Raimbault},
\ead{raimbault@math.jussieu.fr}
\author[samet]{Iddo Samet}
\ead{sameti@math.huji.ac.il}

\address[abert]{Alfr\'ed R\'enyi Institute of Mathematics, POB 127, H-1364
Budapest, Hungary}
\address[bergeron]{Institut de Math\'ematiques de Jussieu
Unit\'e Mixte de Recherche 7586 du CNRS
Universit\'e Pierre et Marie Curie
4, place Jussieu 75252 Paris Cedex 05, France }
\address[biringer]{Yale University
Mathematics Department
PO Box 208283
New Haven, CT 06520}
\address[gelander]{Einstein Institute of Mathematics
Edmond J. Safra Campus, Givat Ram
The Hebrew University of Jerusalem
Jerusalem, 91904, Israel}
\address[nikolov]{Department of Mathematics, Imperial College London, SW7 2AZ, UK}
\address[samet]{Einstein Institute of Mathematics
Edmond J. Safra Campus, Givat Ram
The Hebrew University of Jerusalem
Jerusalem, 91904, Israel}

\medskip
\begin{center}
{\small Received *****; accepted after revision +++++\\
Presented by ¬£¬£¬£¬£¬£}
\end{center}

\begin{abstract}
\selectlanguage{english}
We announce new results concerning the asymptotic behavior of the
Betti numbers of higher rank locally symmetric spaces as their volumes
tend to infinity. Our main theorem is a uniform version of the L\"uck
Approximation Theorem \cite{luck}, which is much stronger than the linear upper
bounds on Betti numbers given by Gromov in \cite{BGS}.

The basic idea is to adapt the theory of local convergence,
originally introduced for sequences of graphs of bounded degree by
Benjamimi and Schramm, to sequences of Riemannian manifolds. Using
rigidity theory we are able to show that when the volume tends to infinity,
the manifolds locally converge to the universal cover in a sufficiently
strong manner that allows us to derive the convergence of the normalized
Betti numbers.
{\it To cite this article: M.
Abert, N. Bergeron, I. Biringer, T. Gelander, N. Nikolov, J. Raimbault, I. Samet, C. R. Acad. Sci. Paris, Ser. I 340 (2011).}

\vskip 0.5\baselineskip

\selectlanguage{francais}
\noindent{\bf R\'esum\'e} \vskip 0.5\baselineskip \noindent
{\bf Comportement des nombres de Betti des espaces localement sym\'etriques.}
Nous annon\c{c}ons de nouveaux r\'esultats concernant le comportement asymptotique des
nombres de Betti des espaces localement sym\'etriques de rang sup\'erieur lorsque leurs volumes
tendent vers l'infini. Notre r\'esultat principal -- une version uniforme du th\'eor\`eme d'approximation
de L\"uck \cite{luck} -- est plus fort que la majoration lin\'eaire en le volume obtenue par Gromov dans
\cite{BGS}.

L'id\'ee de base est d'adapter la th\'eorie de la convergence locale,
initialement introduite pour les suites de graphes de degr\'e born\'e par
Benjamimi et Schramm, \`a des suites de vari\'et\'es riemanniennes. L'utilisation de
th\'eor\`emes de rigidit\'e nous permet de montrer que lorsque le volume tend vers l'infini, les vari\'et\'es convergent localement vers le rev\^etement universel de mani\`ere assez
forte pour en d\'eduire la convergence des
nombres de Betti normalis\'es par le volume.
{\it Pour citer cet article~: M.
Abert, N. Bergeron, I. Biringer, T. Gelander, N. Nikolov, J. Raimbault, I. Samet, C. R. Acad. Sci. Paris, Ser. I 340 (2011).}

\end{abstract}
\end{frontmatter}

\selectlanguage{francais}
\section*{Version fran\c{c}aise abr\'eg\'ee}
Soit $G$ un groupe de Lie simple de rang r\'eel sup\'erieur \`a $2$ et soit $X=G/K$ l'espace sym\'etrique
associ\'e. Une $X$-vari\'et\'e est une vari\'et\'e riemannienne compl\`ete localement isom\'etrique \`a $X$, autrement dit une vari\'et\'e de la forme $%
M=\Gamma \backslash X$, o\`u $\Gamma \leq G$ est un sous-goupe discret sans torsion. On note
$b_k (M)$ le $k$-i\`eme nombre de Betti de $M$ et $\beta_k (X)$ le $k$-i\`eme nombre de Betti $L^2$
de $X$.

\begin{theoreme} \label{TF}
Soit $(M_n)$ une suite de $X$-vari\'et\'es compactes dont le rayon d'injectivit\'e est uniform\'ement
minor\'ee par une constante strictement positive et telle que $\mathrm{vol}(M_{n})\rightarrow
\infty $. Alors on a~:
$$\lim_{n\rightarrow \infty }\frac{b_{k}(M_{n})}{\mathrm{vol}(M_{n})}=\beta
_{k}(X)
$$
pour $0\leq k\leq \dim (X)$.
\end{theoreme}

Margulis conjecture qu'il existe une constante strictement positive qui minore le rayon d'injectivit\'e
de toute $X$-vari\'et\'e compacte (voir \cite[page 322]{margulis:book} et \cite[%
Section 10]{Gel:HV}). Le th\'eor\`eme \ref{TF} s'applique donc conjecturalement \`a toute
suite infinie de $X$-vari\'et\'es compactes distinctes.

Nous montrons d'ailleurs une version faible de la conjecture de Margulis~: pour toute suite $(M_n)$
de $X$-vari\'et\'es compactes distinctes et pour tout r\'eel $r>0$ on a
$$
\lim_{n\rightarrow \infty }\frac{\mathrm{vol}((M_{n})_{<r})}{\mathrm{vol}(M_{n})}=0
$$
o\`u $M_{<r}$ d\'esigne l'ensemble des points de $M$ o\`u le rayon d'injectivit\'e local est
strictement inf\'erieur \`a $r$ (voir Corollary \ref{cor:small-thin} ci-dessous).

On donne une d\'efinition analytique des nombres de Betti $L^2$ de $X$ dans le paragraphe \ref{bettinumbers} ci-dessous. On peut \'egalement les d\'efinir comme suit~: soit $X^*$ le
dual compact de $X$ muni de la m\'etrique Riemannienne induite par la forme de Killing sur $\mathrm{Lie} (G)$. Alors~:
$$
\beta _{k} (X)= \left\{
\begin{array}{ll}
0, & k\neq \frac{1}{2}\dim X \\
\frac{\chi (X^{\ast })}{\mathrm{vol}(X^{\ast })}, & k=\frac{1}{2}\dim X.%
\end{array} \right.
$$
Noter que $\chi (X^{\ast })=0$ sauf si $\mathrm{rank}_{\mathbb{C}}(G)=%
\mathrm{rank}_{\mathbb{C}}(K)$.

L'id\'ee essentielle pour la d\'emonstration du th\'eor\`eme \ref{TF} est d'adapter aux $X$-vari\'et\'es la notion de \og convergence locale \fg \ due \`a Benjamini et Schramm. On montre alors qu'en rang†sup\'erieur toute suite de $X$-vari\'et\'es distinctes de volume fini converge localement vers $X$.
D'autre part, en utilisant la d\'efinition analytique des nombres de Betti on
montre (cette fois sans hypothËse sur le rang de $X$) que la convergence
locale implique la convergence des nombres de Betti normalis\'es par le volume
(sous l'hypoth\`ese d'un rayon d'injectivit\'e uniform\'ement minor\'e).

Le th\'eor\`eme \ref{TF} est faux lorsque le rang r\'eel de $G$ est \'egal \`a $1$. Cependant, dans \cite{paper} nous
d\'emontrerons une version affaiblie du th\'eor\`eme \ref{TF} valable ind\'ependamment du rang r\'eel
de $G$ ainsi qu'une version
forte, \'egalement valable ind\'ependamment du rang r\'eel
de $G$, mais seulement pour les $X$-vari\'et\'es \og de congruence \fg \ associ\'ees \`a une structure
alg\'ebrique donn\'ee de $G$.

\selectlanguage{english}
\section{Introduction}
\label{}

Let $G$ be a simple Lie group with $\mathbb{R}$-rank at least two and let $%
X=G/K$ be the associated symmetric space. An $X$-manifold is a complete
Riemannian manifold locally isomorphic to $X$, i.e. a manifold of the form $%
M=\Gamma \backslash X$ where $\Gamma \leq G$ is a discrete torsion free
subgroup. We denote by $b_{k}(M)$ the $k^{th}$ Betti number of $M$ and
by $\beta _{k} (X)$ the $k^{th}$ $L^{2}$-Betti number of $X$.

\begin{theorem}
\label{thm:main} Let $(M_{n})$ be a sequence of closed $X$-manifolds with
injectivity radius uniformly bounded away from $0$, and $\mathrm{vol}(M_{n})\rightarrow
\infty $. Then:
$$\lim_{n\rightarrow \infty }\frac{b_{k}(M_{n})}{\mathrm{vol}(M_{n})}=\beta
_{k}(X)
$$
for $0\leq k\leq \dim (X)$.
\end{theorem}

Margulis conjectured that there is a uniform lower bound for the injectivity
radius of closed $X$-manifolds (see \cite[page 322]{margulis:book} and \cite[%
Section 10]{Gel:HV}). If this were true, Theorem \ref{thm:main} would apply
to any sequence of distinct, closed $X$-manifolds.

In this direction, our methods lead to the following version of the
Margulis conjecture: for every sequence $(M_{n})$ of distinct, closed $X$%
-manifolds and $r>0$ we have
$$
\lim_{n\rightarrow \infty }\frac{\mathrm{vol}((M_{n})_{<r})}{\mathrm{vol}(M_{n})}=0
$$
where $M_{<r}$ denotes the set of points in $M$ where the local injectivity
radius is less than $r$.

We will define the $L^{2}$-Betti numbers of $X$ analytically in Section \ref%
{bettinumbers} below, but one can also describe them as follows. Let $%
X^{\ast }$ denote the compact dual of $X$ and recall that the Killing form
on $\mathrm{Lie}(G)$ induces a Riemannian structure and a volume form on $%
X^{\ast }$ as well as on $X$. Then
$$
\beta _{k} (X)= \left\{
\begin{array}{ll}
0, & k\neq \frac{1}{2}\dim X \\
\frac{\chi (X^{\ast })}{\mathrm{vol}(X^{\ast })}, & k=\frac{1}{2}\dim X.%
\end{array} \right.
$$
Note also that $\chi (X^{\ast })=0$ unless $\mathrm{rank}_{\mathbb{C}}(G)=%
\mathrm{rank}_{\mathbb{C}}(K)$.

In \cite{paper}, we will obtain analogous results for noncompact $X$%
-manifolds of finite volume and also extend our main theorem to semi-simple
Lie groups. Also, we shall put the main result in the more general context
of counting multiplicities over the unitary dual $\widehat{G}$. Moreover,
for \emph{congruence covers} of a fixed arithmetic $X$-manifold we will
obtain even stronger results on the asymptotic behavior of Betti numbers,
independently of the $\mathbb{R}$-rank of $G$. In particular, Theorem \ref%
{thm:main} holds for any sequence of congruence covers.


It is easy to see that Theorem \ref{thm:main} is false for rank one
symmetric spaces. For instance, suppose that $M$ is a closed hyperbolic $n$%
-manifold such that $\pi _{1}(M)$ surjects on the free group of rank $2$.
Then finite covers of $M$ corresponding to subgroups of $\Z * \Z$ have first
Betti numbers that grow linearly with volume. However, there will be
sublinear growth of Betti numbers in any sequence of covers corresponding to
a chain of finite index normal subgroups of $\pi _{1}(M)$ with trivial
intersection, by \cite{luck}. In \cite{paper}, we shall give a weaker
version of Theorem \ref{thm:main} for rank one spaces.

\section{The main ideas and skeleton of the proof}

\subsection{Local convergence of $G $-spaces}\label{IRSs}

In \cite{localconvergence}, Benjamini and Schramm introduced a notion of probabilistic, or \it local,  \rm convergence for sequences of finite graphs.  Here, we introduce a variant of local convergence for spaces modeled on a fixed Lie group.

Let $G $ be a Lie group and let $\mathrm{Sub}_{G}$  denote the space of closed
subgroups of $G$ equipped with the
Chabauty topology \cite {Chabauty}. As a topological space, $\mathrm{Sub}_{G}$
is compact and $G$ acts on $\mathrm{Sub}_{G}$ continuously by conjugation.



\begin{e-definition}
An \emph{invariant random subgroup} (IRS) of $G$ is a $G$-invariant probability
measure on $\mathrm{Sub}_{G}$.
\end{e-definition}

This notion has been introduced in \cite{agv}. It follows from the compactness of $\mathrm{Sub}_{G} $ that the space of invariant random subgroups of $G$ equipped with the weak* topology is also compact.

The stabilizer of a random point in a continuous probability measure preserving action of
$G$ is an invariant random subgroup.  As an example, if $\Gamma < G $ is a lattice then $G $ acts continuously on the space $\Gamma \backslash G $.  This action preserves a unique probability measure (Haar measure), and we let $\mu_\Gamma $ be the stabilizer of a
random point. Note that the measure $\mu _{\Gamma }$ is supported
on the conjugacy class of $\Gamma $.  We also let $\mu _{\mathrm{id}}$ be the
measure supported on the trivial subgroup of $G$ and $\mu _{G}$
be the measure supported on the element $G\in \mathrm{Sub}_{G} $.

Using Borel's density theorem \cite {borel}, one can prove that if $G $ is simple then every IRS other than $\mu_G $ is supported on discrete subgroups of $G $.  It is well-known that for discrete subgroups $\Lambda  < G $, convergence in the Chabauty topology is equivalent to Gromov Hausdorff convergence of the quotient manifolds $\Lambda \backslash G $.  Convergence of invariant random subgroups then has the following interpretation, which is an exact analogue of Benjamini-Schramm convergence of finite graphs.

\begin{lemma}[Local convergence] \label{lem:loccv}
Suppose that $\mu_1,\mu_2,\ldots,\mu_\infty$ are IRSs supported on discrete subgroups of a Lie group $G$. Then the following are equivalent:
\begin{enumerate}
\item $\mu_n$ converges weakly to $\mu_\infty$.
\item For every pointed Riemannian manifold $(M,p)$ and $R >0 $, and arbitrarily small $\eps,\delta > 0$, the probability that for a $\mu_n$-random $\Gamma \in \mathrm{Sub}_G$ the pointed ball $(B_{\Gamma \backslash G}([\mathrm{id}],R),[\mathrm{id}])$ is $(1+\delta,\eps)$-quasi-isometric to $(B_M(p,R),p)$ converges to the $\mu_\infty$-probability of this event.
\end{enumerate}
\end{lemma}

\subsection {Local convergence in higher rank}  Suppose now that $G $ is a simple Lie group with $\R $-rank at least 2 and associated symmetric space $X = G/K $.

In this case, one can classify all the ergodic invariant random subgroups of $G $.

\begin{e-proposition}\label{weak}
Every ergodic invariant random subgroup of $G$ is equal to either $\mu _{\mathrm{id}}$, $\mu _{G}$ or to $\mu _{\Gamma }$ for some
lattice $\Gamma $ in $G$.
\end{e-proposition}

Proposition \ref{weak} is implicitely contained in \cite[\S 3]{stuckzimmer}. The proof relies on the rigidity theory of higher rank lattices: in particular, on the St\"{u}ck--Zimmer theorem \cite[Theorem 2.1]{stuckzimmer}, Borel's Density theorem \cite{borel} and the Margulis normal subgroup theorem \cite[Chapter VIII]{margulis:book}.

Using Proposition \ref {weak} and a quantitative version of the Howe--Moore theorem due to Howe and Oh \cite{Oh}, one can prove the following key result:

\begin{theorem}\label {weak2} $\mu _{\mathrm{id}}$ is the only accumulation point of
$\left\{ \mu _{\Gamma }\mid \Gamma \mbox{ is a lattice in }G\right\}.$
\end{theorem}


The geometric interpretation of weak convergence given in Theorem \ref{lem:loccv} implies that if $\Gamma_n < G $ is a sequence of lattices such that $\mu_{\Gamma_n} \to \mu_{\mathrm{id}} $, then the injectivity radius of $\Gamma_n \backslash X$ at a random point goes to infinity asymptotically almost surely.  Therefore, Theorem \ref {weak2} has the following corollary.

\begin{corollary}\label{cor:small-thin}
Suppose that $G$ is a simple Lie group with $\mathbb{R}$-rank at least $2$ and associated symmetric space $X = G/K $. If $M_n$ is a sequence of distinct, finite volume $X $-manifolds, then for every $r>0$ we have
$$
 \lim_{n\to\infty}\frac{\vol ((M_n)_{<r})}{\vol(M_n)}=0,
$$
where $M_{<r}:=\{x\in M \ | \ \mathrm{InjRad}_M(x)<r\}$ is the $r $-thin part of $M $.
\end{corollary}



\subsection{The Laplacian and heat kernel on differential forms}


As before, let $X = G/K $ be the symmetric space associated to a simple Lie group.  Denote by $e^{-t \Delta_k^{(2)}}
\in \mathrm{End}(\Omega^k_{(2)}(X))$
the bounded operator  (cf. \cite{BarbaschMoscovici}) on $L ^ 2 $ $k $-forms defined by
the fundamental solution of the heat equation:
$$\left\{
\begin{array}{l}
\Delta_k^{(2)} P_t = -\frac{\partial}{\partial t} P_t , \ \ \ t>0 \\
P_0 = \delta
\end{array} \right.$$
where $\delta$ is the Dirac distribution.
It is an integral operator with kernel $e^{-t\Delta_k^{(2)}}(x, y)$ (the heat kernel):  that is, $(e^{-t\Delta_k^{(2)}} f) (x) = \int_{X} e^{-t\Delta_k^{(2)}} (x,y) f(y) \, dy, \ \ \forall f \in  \Omega^k_{(2)}(X).$

From \cite[Lemma 3.8]{BV} we deduce:

\begin{corollary}\label{cor:heat}
Let $m >0$. There exists a positive constant $c_1=c_1(G,m)$ such that
$$|| e^{-t \Delta_k^{(2)}} (x,y) || \leq c_1 e^{- d (x, y)^2 / c_1 } ,  \quad  |t| \leq m.$$
\end{corollary}

Now let $M=\Gamma\backslash X$ be a compact  $X$-manifold. Let $\Delta_k$ be the Laplacian on differentiable $k$-forms on
$M$. It is a symmetric, positive definite, elliptic operator with pure point spectrum.
Write $e^{-t\Delta_k} (x,y)$ ($x,y \in M$) for the integral kernel of the heat equation of $k$-forms on $M$, then for each positive $t$ we have:
\begin{equation} \label{eq:heatkernels}
e^{-t\Delta_k} (x,y) = \sum_{\gamma \in \Gamma} (\gamma_y)^* e^{-t\Delta_k^{(2)}} (\widetilde{x} , \gamma \widetilde{y}),
\end{equation}
where $\widetilde{x}, \widetilde{y}$ are lifts of $x,y$ to $X$ and by $(\gamma_y)^*$,
we mean pullback by the map $(x,y) \mapsto (x, \gamma y)$.
The sum converges absolutely and uniformly for $\tilde{x}, \tilde{y}$ in compacta; this follows from Corollary \ref{cor:heat}, together with the following estimate:
\begin {equation}\label {eq:numtranslates}
 \left|\{ \gamma \in \Gamma \; : \; d(\tilde x,\gamma \tilde y) \leq r \} \right| \leq c_2 e ^ {c_2 r} \,\mathrm{InjRad}_M(x)^{-d} .
\end {equation}
Here, $c_2=c_2(G)$ is some positive constant and $d=\dim(X)$.




\subsection{($L^2$-)Betti numbers}
\label {bettinumbers}
The trace of the heat kernel $e^{-t \Delta_k^{(2)}}(x,x)$ on the diagonal
is independent of $x \in X$, being $G$-invariant. We denote it by $$\mathrm{Tr} \,e^{-t \Delta_k^{(2)}}:=\mathrm {tr} \, e^{-t \Delta_k^{(2)}}(x,x)$$
and set
$$\beta_k (X) :=  \lim_{t \rightarrow \infty} \mathrm{Tr} e^{-t \Delta_k^{(2)}} .$$
Recall also that the usual Betti numbers of $M$ are given by
$$b_k (M) = \lim_{t \rightarrow \infty} \int_M \mathrm{tr} \ e^{-t \Delta_k}(x,x)dx.$$


The following is a consequence of Corollary \ref {cor:heat} and Equations (\ref {eq:heatkernels}) and (\ref {eq:numtranslates}).

\begin{lemma} \label{lem:hk}
Let $m >0$ be a real number. There exists a constant $c = c(m,G)$ such that for any $x \in M$ and $t\in (0,m]$,
$$
 \left| \mathrm{tr} \ e^{-t\Delta_k } (x,x) - \mathrm{Tr} e^{-t \Delta_k^{(2)}} \right| \leq c \cdot \mathrm{InjRad}_M(x)^{-d}.
$$
\end{lemma}

\subsection{Concluding the proof of Theorem \ref{thm:main}}\label{par:hyp}
Let now $M_n=\Gamma_n\backslash X$ as in the statement of \ref{thm:main}. Since the injectivity radius is uniformly bounded away from $0$
it follows from Corollary \ref{cor:small-thin} that
$\frac{1}{\mathrm{vol} (M_n)} \int_{(M_n)_{< r}} \mathrm{InjRad}_{M_n}(x)^{- d} dx \rightarrow 0
$ for every $r>0$.
Hence by Lemma \ref{lem:hk}:

$$\frac{1}{\vol(M_n)} \int_{M_n} \mathrm{tr} \ e^{-t\Delta_k^{M_n} } (x,x) \, dx \, \longrightarrow \, \mathrm{Tr} e^{-t \Delta_k^{(2)}}$$
uniformly for $t$ on compact subintervals of $(0, \infty)$. Since each function in the limit above is decreasing as a function of $t$, we deduce that

$$
\limsup_{n \rightarrow \infty} \frac{b_k (M_n)}{\mathrm{vol} (M_n)} \leq \beta_k (X).
$$
Finally, using that the usual Euler characteristic is equal to its $L^2$ analogue and that $\Delta_k^{(2)}$ has
zero kernel if $k \neq \frac12 \dim X$ we derive Theorem \ref{thm:main}.

\section*{Acknowledgements}
M.A. has been supported by the grant IEF-235545.
N.B. is a member of the Institut Universitaire de France.  I. B. is partially supported by NSF postdoctoral fellowship DMS-0902991. T.G. acknowledges support of the European Research Council (ERC)/ grant agreement 203418 and the ISF. I.S. acknowledges the support of the ERC of T. Gelander
/ grant agreement 203418.

\def\cprime{$'$} \def\cprime{$'$} \def\cprime{$'$}

\end{document}